\documentclass[12pt]{article}

\usepackage{amssymb}

\newtheorem {theorem} {Theorem}%[section]

\newtheorem {proposition} [theorem]{Proposition}

\def\div{{\rm div}}
\def\cn{{\rm cn}}
\def\sn{{\rm sn}}
\def\dn{{\rm dn}}
\def\am{{\rm am}}

\title{On the stability of periodic orbits for \\ differential systems
in $\mathbb{R}^n$.\thanks{The first author is partially supported
by a DGICYT grant number MTM2005-06098-C02-01 and by a CICYT grant
number 2005SGR 00550. The second and third authors are partially
supported by a DGICYT grant number MTM2005-06098-C02-02. The first
and second authors are also supported by the CRM Research Program:
On Hilbert's 16th Problem.}}

\author{{\sc Armengol Gasull$^{\ (1)}$, H\'ector Giacomini$^{\ (2)}$} \\ {\sc and Maite Grau$^{\ (3)}$}}
\date{}

\begin{document}
\maketitle \vspace{-1.0cm}

\begin{center}
\noindent {\normalsize{$^{\ (1)}$ Departament de Matem\`atiques.
Universitat Aut\`onoma de Barcelona. \\ 08193 - Bellaterra
Barcelona, SPAIN. \\ {\rm E--mail:} {\tt gasull@mat.uab.es}
\\ $ $ \\
$^{\ (2)}$ Lab. de Math\'ematiques et Physique Th\'eorique. CNRS
UMR 6083. \\ Facult\'e des Sciences et Techniques. Universit\'e de
Tours. \\ Parc de Grandmont, 37200 Tours, FRANCE.
\\ {\rm E-mail:} {\tt giacomini@phys.univ-tours.fr}
\\ $ $ \\
$^{\ (3)}$ Departament de Matem\`atica. Universitat de Lleida. \\
Avda. Jaume II, 69. 25001 Lleida, SPAIN. \\ {\rm E--mail:} {\tt
mtgrau@matematica.udl.es}}}
\end{center}
\begin{abstract} We consider an autonomous differential system in
$\mathbb{R}^n$ with a periodic orbit and we give a new method for
computing the characteristic multipliers associated to it. Our
method works when the periodic orbit is given by the transversal
intersection of $n-1$ codimension one hypersurfaces and is an
alternative to the use of the first order variational equations.
We apply it to study the stability of the periodic orbits in
several examples, including a periodic solution found by Steklov
studying the rigid body dynamics.
\end{abstract}

{\small{\noindent 2000 {\it AMS Subject Classification:}  34D08; 37D05, 70E50. \\
\noindent {\it Key words and phrases:} Periodic orbit, characteristic
multipliers, invariant curve, rigid body dynamics, Mathieu's equation, Steklov
periodic orbit.}}

\section{Introduction and statement of the results. \label{sect1}}

Consider a differential system in $\mathbb{R}^n$, with $n \geq 2$, given by:
\begin{equation}
\frac{d \mathbf{x}}{dt}\ = \ \mathbf{X}(\mathbf{x}), \label{eq1}
\end{equation}
where $\mathbf{X}: \mathcal{U} \subseteq \mathbb{R}^n \to \mathbb{R}^n$ is a
$\mathcal{C}^1$ function in some non--null open set $\mathcal{U} \subseteq
\mathbb{R}^n$ and $t$ is a real independent variable. We assume that system
(\ref{eq1}) exhibits a periodic orbit $\Gamma:=\{ \gamma(t)\, | \, 0 \leq t <T
\} \subset \mathcal{U}$ with period $T>0$. As usual we will  denote by
$\phi(t,p)$ the flow solution of (\ref{eq1}) such that $\phi(0,p)=p,$  by $x_i$
the $i^{th}$ component of the point $\mathbf{x}$, that is, $\mathbf{x}
=(x_1,x_2, \ldots, x_n)^{\rm T}$ and analogously $\mathbf{X}=(X_1,X_2, \ldots,
X_n)^{\rm T}$, where ${\rm T}$ denotes transposition.
\newline

It is well known that to determine the behavior of the flow near
$\Gamma$, a first step is to get the characteristic multipliers
associated to this orbit. These multipliers are usually obtained
through the study of the first order variational equations. In
this paper we propose an alternative way for obtaining them which
works when $\Gamma$ is given by the transversal intersection of
$n-1$ codimension one hypersurfaces.
\newline

Recall that the behavior near $\Gamma$ is given by the Poincar\'e map, which is
defined in a section $\Sigma$ through a point $p \in \Gamma,$ where $\Sigma$
passes through $p$ and is a local smooth manifold of dimension $n-1$
transversal to $\Gamma.$ Given a section $\Sigma$,  the {\em Poincar\'e map} is
defined as:
\[\begin{array}{llll} \Pi: & \Sigma  & \longrightarrow & \Sigma \\ & q
& \mapsto & \phi(\tau(q),q), \end{array} \] where $\tau(q)$ is the unique real
function such that $\phi(\tau(q),q) \in \Sigma$ and $\displaystyle \lim_{q \to
p} \tau(q) = T$. \newline

 The $n-1$ eigenvalues of $D\Pi(p)$ are independent of
$p$ and $\Sigma$ and are  called the {\em characteristic multipliers} of
$\Gamma.$ It is well known that the stable (respectively, unstable) manifold
associated to $\Gamma$ has dimension the number of characteristic multipliers
with modulus smaller than 1  (respectively, bigger that 1). It is also well
known, see \cite{Liapunov}, that if all the characteristic multipliers have
modulus lower than or equal to $1$, then $\Gamma$ is Liapunov stable and  if
all the characteristic multipliers have modulus strictly lower than $1$, then
it is asymptotically stable. \newline

Our  result is  motivated from a previous result given in
\cite{hyper}  for planar differential systems.  Recall that in the
planar case there is only one characteristic multiplier which is
given by:
\[ \Pi'(p) \, = \, \exp \left\{ \int_{0}^{T} \div(\mathbf{X})(\gamma(t)) \, dt
\right\}, \] where $\div(\mathbf{X})={\partial X_1}/{\partial x_1}  + {\partial
X_2}/{\partial x_2}$ is the {\em divergence} of the system, see for instance
\cite[p. 214]{Perko}.  In  \cite{hyper}, the authors give an alternative
expression for the same value $\Pi'(p)$. Let us consider an {\em invariant
curve} $f(x_1,x_2)=0$ for a planar system (\ref{eq1}), that is a curve defined
by a real $\mathcal{C}^1(\mathcal{U})$--function $f(x_1,x_2)$, for which there
exists a function $k(x_1,x_2): \mathcal{U} \subseteq \mathbb{R}^2 \to
\mathbb{R}$ of class $\mathcal{C}^1$ satisfying that: \[
Df(\mathbf{x})\,\mathbf{X}(\mathbf{x})=\nabla f (\mathbf{x})\cdot
\mathbf{X}(\mathbf{x})^{\rm{T}} \, = \, k (\mathbf{x}) \, f (\mathbf{x}),
\] where $\nabla f \, = \, \left( \partial f / \partial x_1,
\partial f / \partial x_2 \right)$. Under these assumptions, the
function $k(x_1,x_2)$ is called the {\em cofactor} associated to
the invariant curve $f(x_1,x_2)=0$.
\begin{proposition} {\sc \cite{hyper}} Consider system {\rm (\ref{eq1})} with $n=2$ and a
periodic orbit $\Gamma=\{ \gamma(t)\, | \, 0 \leq t <T \}$ contained in a
planar invariant curve $f(x_1,x_2)=0$. If $\nabla f$ does not vanish on
$\Gamma$, then
\[ \Pi'(p) \, = \, \exp \left\{ \int_{0}^{T} k(\gamma(t)) \, dt
\right\}. \] \label{prophyper} \end{proposition} The above result allows to
compute the characteristic multiplier associated to $\Gamma$ through the
integration over the cofactor associated to an invariant curve containing the
periodic orbit. In Theorem \ref{th1}, we  extend this result to  periodic
orbits in $\mathbb{R}^n$ given as the transversal intersection of $n-1$
codimension one hypersurfaces.\newline

\begin{theorem} Let $\Gamma=\{ \gamma(t)\, | \, 0 \leq t <T
\}$ be  a $T-$periodic solution of {\rm (\ref{eq1})}. Consider a
smooth function $\mathbf{f}: \mathcal{U} \subseteq \mathbb{R}^n
\to \mathbb{R}^{n-1}$, $\mathbf{f}= \left( f_1, f_2, \ldots,
f_{n-1} \right)^{\rm T}$, such that: \vspace{0.2cm}
\par \vspace{0.2cm} - $\Gamma$ is contained in $\displaystyle
\bigcap_{i=1,2, \ldots, n-1} \{f_i(\mathbf{x})=0\}$,
\par \vspace{0.2cm}
- the crossings of all the manifolds $\{f_i(\mathbf{x})=0\}$ for $i=1,2,\ldots,
n-1$ are transversal over $\Gamma$,
\par \vspace{0.2cm}
- there exists a $(n-1) \times (n-1)$ matrix
$\mathbf{k}(\mathbf{x})$ of real functions satisfying:
\begin{equation}
 D\mathbf{f}(\mathbf{x}) \, \mathbf{X}(\mathbf{x}) = \mathbf{k}(\mathbf{x}) \,
\mathbf{f}(\mathbf{x}). \label{eq2}
\end{equation}
\par \vspace{0.3cm}
Let $\mathbf{v}(t)$ be the $(n-1) \times (n-1)$ fundamental matrix solution of
\begin{equation} \frac{d \mathbf{v}(t)}{dt} \, = \, \mathbf{k}(\gamma(t)) \,
\mathbf{v}(t),\qquad \mathbf{v}(0) = {\rm Id}. \label{eq3}
\end{equation} Then  the characteristic
multipliers of $\Gamma$ are the eigenvalues of $\mathbf{v}(T)$. \label{th1}
\end{theorem}

Indeed in the proof of this theorem, we will deduce that the
matrices $\mathbf{v}(T)$ and $D\Pi(p)$ are conjugated. On the
other hand, recall that the usual way of computing the
characteristic multipliers of $\Gamma$ lies in finding the
fundamental matrix $u(T),$ solution of the variational equation:
\begin{equation} \frac{d
\mathbf{u}(t)}{dt} \, = \, D\mathbf{X}(\gamma(t)) \, \mathbf{u}(t), \qquad
\mathbf{u}(0)={\rm Id}, \label{eq33}
\end{equation}
which is a $n\times n$ matrix having eigenvalues 1 and  the $n-1$
characteristic multipliers of $\Gamma.$
\newline

Note that the assumptions stated in our theorem are
straightforward generalizations of the notion of invariant curve
containing the periodic orbit $\Gamma$ for a system in dimension
$n$ and the integration of the cofactor, which in $\mathbb{R}^n$
consists on solving the linear differential system $\mathbf{v}'(t)
\, = \, \mathbf{k}(\gamma(t)) \, \mathbf{v}(t).$
\newline

Next section is devoted to prove the above theorem.  Finally, in
Section \ref{sect2} we apply it to study the stability of the
periodic orbits of several differential equations.  Example 1 is
devoted to study a 4-dimensional polynomial system, that  includes
two systems studied in \cite{dekleine}, exhibiting an explicit
periodic orbit.  Example 2 shows a 3-dimensional polynomial system
for which we prove that the stability of the given periodic orbit
is equivalent to the study of the stability of the Mathieu's
equation. The last example deals with a more involved case, a
6-dimensional system that controls the dynamics of a rigid body.
In our approach, and by using the existence of three independent
first integrals, we  reduce the study of the stability of the
Steklov periodic orbit to the study of a second order linear
differential equation. This Steklov orbit is introduced in
\cite{Steklov}. It is worth to say that in all the above examples
we have also tried to use the usual approach, namely the
variational equations, for studying the stability of the given
periodic orbit and we have found that the result of Theorem
\ref{th1} makes the computations easier.

\section{Proof of Theorem {\rm \ref{th1}}\label{sect15}}

Fix a point $ p \in\Gamma.$  To prove the result we consider the Poincar\'{e}
section given by the  $(n-1)-$dimensional orthogonal hyperplane to $\mathbf{X}(
p ),$
$$\Sigma=<{\mathbf {X}}( p)>^{\perp}=<{\bf e}_1,\ldots{\bf e}_{n-1}>.$$

Take a new system of coordinates centered at $ p $ with basis ${ \bf
e}_1,\ldots{\bf e}_{n-1}, { \bf X}({ p}).$ These new coordinates write as ${\bf
y}= A({\bf x}-{p}),$  for some invertible matrix $A.$ By using them the
differential equation (\ref{eq1}) is converted into
\begin{equation}
\frac{d \mathbf{y}}{dt}\ = \ \mathbf{Y}(\mathbf{y}):=\
A\mathbf{X}(A^{-1}\mathbf{y}+p), \label{eq1new}
\end{equation}
and the manifold containing $\Gamma$ is given by
$\mathbf{g}(\mathbf{y}):=\mathbf{f}(A^{-1}\mathbf{y}+p)=\mathbf{0},$
which satisfies
\begin{equation}
 D\mathbf{g}(\mathbf{y}) \, \mathbf{Y}(\mathbf{y}) = \mathbf{\tilde k}(\mathbf{y}) \,
\mathbf{g}(\mathbf{y}), \label{gnew}
\end{equation}
where $\mathbf{\tilde k}(\mathbf{y}):=\mathbf{ k}(A^{-1}\mathbf{y}+p ),$
because by using (\ref{eq2}), the following equality holds:
$$
 D\mathbf{f}(A^{-1}\mathbf{y}+p) \,A^{-1}\,A\, \mathbf{X}(A^{-1}\mathbf{y}+p) = \mathbf{ k}(A^{-1}\mathbf{y}+p) \,
\mathbf{f}(A^{-1}\mathbf{y}+p).
$$

Recall that  $\phi(t, p )$ is the solution of (\ref{eq1}) such
that when $t=0$ passes through $ p.$ Let $\psi(t, q )$ be the
solution of (\ref{eq1new}) such that when $t=0$ passes through $
q.$ Then $\psi(t, q )=A[\phi(t,A^{-1} q +p)- p ].$

In these new coordinates note that $\Sigma=\{{y_n=0}\}$ and if $ q
=(y_1,y_2,\ldots,y_{n-1},0)^{\rm T}\in\Sigma$ then the Poincar\'{e} map
$\Pi:\Sigma\to\Sigma$ writes as
$$\Pi( q )\, = \, \left( \psi_1(\tau( q ), q ),\,
\psi_2(\tau( q ), q ), \, \ldots, \, \psi_{n-1}(\tau( q ), q ),\,
0 \right)^{\rm{T}},$$ where $\tau( q )$ is precisely  the time
such that $\psi_n((\tau( q ), q ))=0,$ which is known to be a
smooth function. Indeed we can identify $\Pi$ with a map $\tilde
\Pi$ from $R^{n-1}$ into itself with variables $y_1,\ldots,
y_{n-1},$ defined as $$\tilde \Pi(y_1,\ldots,y_{n-1})=\left(
\Pi_1(y_1,\ldots,y_{n-1},0),\, \ldots, \,
\Pi_{n-1}(y_1,\ldots,y_{n-1},0) \right)^{\rm{T}}.$$

Observe also that the hypotheses on the transversal cutting of
$\Gamma$ and the hypersurfaces $f_i(\mathbf{ x})=0$ imply that
extended map, $\mathbf{ \hat f}=(f_1,f_2,\ldots,f_{n-1},f_n),$
where $f_n(\mathbf{x}):=X_1(p)x_1+X_2(p)x_2+\cdots X_n(p)x_n$ is
such that $\det(D\mathbf{\hat f}(p))\ne0.$ This information
translated to the function $\mathbf{g}(\mathbf{y})$ implies that
the matrix $D\mathbf{g}({0})$ has rank $n-1$ and that a minor with
determinant different from zero is the one corresponding to the
partial derivatives with respect to $y_1,y_2,\ldots, y_{n-1}.$

On the other hand,  by using (\ref{gnew}) we have that:
\begin{eqnarray*} \displaystyle \frac{\partial\mathbf{g} \left(
\psi(t,q) \right)}{\partial t} & = & \displaystyle D\mathbf{g}
\left( \psi(t,q) \right) \, \frac{\partial \left( \psi(t,q)
\right)}{\partial t} \, = \, D\mathbf{g} \left( \psi(t,q) \right)
\, \mathbf{Y} \left( \psi(t,q) \right) \vspace{0.2cm}
\\ & = & \displaystyle \mathbf{\tilde k}\left(
\psi(t,q) \right) \, \mathbf{g} \left( \psi(t,q) \right),
\end{eqnarray*} for any point $q$ in the domain of definition of (\ref{eq1new}).

Let $\mathbf{v}(t;q)$ be the fundamental matrix solution of
\begin{equation}\label{ult}
 \frac{d \mathbf{v}(t)}{dt} \, = \, \mathbf{\tilde k}\left( \psi(t,q)
\right) \, \mathbf{v}(t).
\end{equation} Thus $ \mathbf{g} \left( \psi(t,q) \right)
\, = \, \mathbf{v}(t;q) \, \mathbf{g}(q)$, because the function
$\mathbf{g} \left( \psi(t,q) \right)$ is a solution of the same
linear differential system and satisfies that $\mathbf{g} \left(
\psi(0,q) \right)=\mathbf{g} \left(q \right).$
\newline

Consider a point $q \in \Sigma.$  In these coordinates, we have
$q=(q_1,q_2,\ldots,q_{n-1},0)^{\rm T}$. We define
$z=(q_1,q_2,\ldots,q_{n-1})^{\rm T}$  and $h(z)=(z^{\rm T},0)^{\rm
T}=q$ like the inclusion of $z$ as a point in the domain of
definition of (\ref{eq1new}). Since $\Pi(q) = \psi(\tau(q),q),$ by
the above result we have that $ \mathbf{g} \left( \Pi(q) \right) =
\mathbf{v} \left(\tau(q),q \right) \, \mathbf{g}(z),$ or
equivalently that
$$\mathbf{g}\left( h(\tilde\Pi(z)) \right) = \mathbf{v} \left(\tau(h(z));h(z) \right)
\, \mathbf{g}(h(z)).$$

We differentiate the previous identity with respect to $z$:
\begin{eqnarray*}
 D\mathbf{g}\left( h(\tilde\Pi(z)) \right) \, Dh(\tilde\Pi(z)) \, D\tilde\Pi(z) &=& D\left[\mathbf{v}
\left(\tau(h(z));h(z) \right)\right] \, \mathbf{g}(h(z)) \, +
\\ &+& \, \mathbf{v} \left(\tau(h(z));h(z)\right) \, D\mathbf{g}(h(z))\, Dh(z).
\end{eqnarray*}
By evaluating at $z=\mathbf{0}$, which corresponds exactly to the
point $p\in \Gamma$, and using that
$\mathbf{g}(h(\mathbf{0}))=\mathbf{0}$,
$h(\mathbf{0})=(\mathbf{0}^{\rm T},0)^{\rm T},$
$\tau(\mathbf{0})=T,$ $\tilde\Pi(\mathbf{0})=\mathbf{0}$ and the
expression of $h$, we obtain that:
\[ D_z\mathbf{g}(\mathbf{0}) \, D\tilde\Pi(\mathbf{0}) \, = \,  \mathbf{v}
\left(T;h(\mathbf{0})
\right) \, D_z\mathbf{g}(\mathbf{0}), \] where
$D_z\mathbf{g}(\mathbf{0}):= D\mathbf{g}(h(\mathbf{0})) \,
Dh(\mathbf{0})$ is precisely the invertible squared matrix formed
by the derivative of $\mathbf{g}$ only with respect to
$q_1,\ldots,q_{n-1}.$  Hence we have that the characteristic
multipliers associated to $\Gamma$ coincide with the eigenvalues
of $\mathbf{v}(T;\mathbf{0})$ because the matrices
$D\tilde\Pi(\mathbf{0})$ and $\mathbf{v}(T;h(\mathbf{0}))$ are
similar.

Finally note that $\mathbf{\tilde k}(
\psi(t,h(\mathbf{0})))=\mathbf{ k}(A^{-1}
\psi(t,h(\mathbf{0}))+p)=\mathbf{ k}( \phi(t,p)).$ Thus equation
(\ref{ult}) is the same that
$$
 \frac{d \mathbf{v}(t)}{dt} \, = \,\mathbf{ k}( \gamma(t)) \, \mathbf{v}(t),\label{var1}
$$
and $\mathbf{v}\left(T;h(\mathbf{0}) \right)$ coincides with the
matrix $\mathbf{v}(T)$ as defined in (\ref{eq3}), as we wanted to
prove.

\section{Some examples and applications. \label{sect2}}

$\quad \,$ {\bf Example 1.} The following example is an extension  of two
systems  extracted from \cite{dekleine}. The goal of the examples given in
\cite{dekleine} is to illustrate that the asymptotic stability of a periodic
orbit of a system (\ref{eq1}) is not determined by the eigenvalues of the
matrix defining the first variational equation. \par We consider the following
differential system in $\mathbb{R}^4$:
\begin{equation}
\begin{array}{lll}
\displaystyle \dot{x} & = & \displaystyle - y - x \left( x^2 + y^2
-1 \right), \vspace{0.2cm} \\
\displaystyle \dot{y} & = & \displaystyle  x - y \left( x^2 + y^2
-1 \right), \vspace{0.2cm} \\
\displaystyle \dot{z} & = & \displaystyle -w - s \, z \left( z^2 +
w^2 -1 \right) - s \, k \left(x-z \right) , \vspace{0.2cm} \\
\displaystyle \dot{w} & = & \displaystyle z - s \, w \left( z^2 +
w^2 -1 \right) - s \, k \left(y-w \right),
\end{array} \label{eq5}
\end{equation}
with $s, k$ real parameters. This system coincides with the first example given
in \cite{dekleine} when $s=1$ and with the second example when $s=-1$. We note
that this system always exhibits the periodic orbit $\Gamma \, = \, \left\{
\gamma(t)  \, : \, 0 \leq t < 2 \pi \right\}, $ where $\gamma(t) \, = \, \left(
\cos(t), \sin(t), \cos(t), \sin(t) \right)$. We are able to compute all the
characteristic multipliers associated to $\Gamma$ for any real value of the
parameters $s$ and $k$. We consider the hypersurfaces $f_i(x,y,z,w)=0$,
$i=1,2,3$, given by: $f_1(x,y,z,w) \, = \, x^2 + y^2 -1$, $f_2(x,y,z,w) \, = \,
x - z$ and $f_3(x,y,z,w) \, = \, y - w$. We denote by $\mathbf{f} \, = \,
\left( f_1, f_2, f_3 \right)$. It is easy to see that $\Gamma$ is contained in
the intersection of these three hypersurfaces and that the crossings of these
hypersurfaces are normal over the periodic orbit $\Gamma$, as the computation
of the following determinant shows:
\[ \left| \begin{array}{c} \nabla f_1(\gamma(t)) \vspace{0.2cm} \\
\nabla f_2(\gamma(t)) \vspace{0.2cm} \\ \nabla f_3(\gamma(t)) \vspace{0.2cm} \\
\mathbf{X}(\gamma(t)) \end{array} \right| \, = \, \left| \begin{array}{cccc} 2
\cos(t) & 2 \sin(t) & 0 & 0 \vspace{0.2cm}
\\ 1 & 0 & -1 & 0  \vspace{0.2cm} \\ 0 & 1 & 0 & -1  \vspace{0.2cm}
\\ - \sin(t) & \cos(t) & - \sin(t) & \cos(t) \end{array} \right| \, =
\, 4 \, \neq \, 0.  \] Straightforward computations show that $
D\mathbf{f}(\mathbf{x}) \, \mathbf{X}(\mathbf{x}) = \mathbf{k}(\mathbf{x}) \,
\mathbf{f}(\mathbf{x}),$ with the following matrix of cofactors:
\[ \mathbf{k}(x,y,z,w):= \left( \begin{array}{ccc} -2 \left( x^2 +
y^2 \right) & 0 & 0 \vspace{0.2cm} \\ s \, z - x &  s \, k  - s \,
z \left( x + z \right) & -1 - s \, z \left( y + w \right)
\vspace{0.2cm} \\ s \, w - y & 1 - s \, w \left( x + z \right) & s
\, k - s \, w \left( y + w \right) \end{array} \right), \] which
evaluated on the periodic orbit $\Gamma$ reads for: \[
\mathbf{k}(\gamma(t)) \, = \, \left( \begin{array}{ccc} -2 & 0 & 0
\vspace{0.2cm} \\ \left( s-1 \right) \cos(t) &  s \left( k - 2
\cos^2(t) \right) & -1-2 s \cos(t) \sin(t)  \vspace{0.2cm} \\
\left( s-1 \right) \sin(t) & 1 - 2 s \cos(t) \sin(t) & s  \left( k - 2
\sin^2(t) \right) \end{array} \right). \] The fundamental matrix solution of
the linear equation $\mathbf{v}'(t) \, = \, \mathbf{k}(\gamma(t)) \,
\mathbf{v}(t)$ is:
\begin{itemize}
\item when $2 - 2 s + s k \neq 0$,
\[ \mathbf{v}(t) \, = \, \left( \begin{array}{ccc} \displaystyle
e^{-2 t} & \displaystyle 0 & \displaystyle 0 \vspace{0.2cm} \\
\displaystyle \left( s -1 \right) e^{-2 t} \left( \frac{ e^{(2 - 2
s + k s) \, t} -1}{2 - 2 s + k s} \right) \cos(t)
& \displaystyle e^{(k-2) s t} \cos(t) & \displaystyle - e^{k s t} \sin(t) \vspace{0.2cm} \\
\displaystyle \left( s -1 \right) e^{-2 t} \left( \frac{ e^{(2 - 2
s + k s) \, t} -1}{2 - 2 s + k s} \right) \sin(t) & \displaystyle
e^{(k-2) s t} \sin(t) & \displaystyle  e^{k s t} \cos(t)
\end{array} \right), \]
\item when $ 2 - 2 s + s k = 0$, we put $k=2(s-1)/s$:
\[ \mathbf{v}(t) \, = \, \left( \begin{array}{ccc} \displaystyle
e^{-2 t} & \displaystyle 0 & \displaystyle 0 \vspace{0.2cm} \\
\displaystyle \left( s -1 \right) e^{-2 t}\, t \, \cos(t)
& \displaystyle e^{-2 t} \cos(t) & \displaystyle - e^{2 ( s -1) t} \sin(t) \vspace{0.2cm} \\
\displaystyle \left( s -1 \right) e^{-2 t} \, t \, \sin(t) &
\displaystyle e^{-2 t} \sin(t) & \displaystyle  e^{2(s-1) t}
\cos(t) \end{array} \right). \]
\end{itemize}
We can compute the eigenvalues of the matrix $\mathbf{v}(2 \pi)$,
which by Theorem \ref{th1}, correspond to the characteristic
multipliers associated to $\Gamma$. These eigenvalues are:
\begin{itemize}
\item when $2 - 2 s + s k \neq 0$: $e^{- 4 \pi}$, $e^{2 s (k-2)
\pi}$ and $e^{2 k s \pi}$,
\item when $k=2(s-1)/s$: $e^{- 4 \pi}$, $e^{- 4 \pi}$ and $e^{4
\pi(s-1)}$.
\end{itemize}
Therefore, for instance in the first case ($2 - 2 s + s k \neq 0$) when $s \, k
>0$ or $s \, (k-2) >0$ we have that $\Gamma$ is unstable and, in
the second case ($2 - 2 s + s k = 0$) when $s<1$ then $\Gamma$ is Liapunov
unstable.
\newline

{\bf Example 2.} We give an example related to the Mathieu's
equation. The Mathieu's equation is a particular case of the
Hill's equation and it has the form: \[ v''(t) \, + \, \left(a + 2
q \cos (2 t) \right) v(t) \, = \, 0, \] where $a,q$ are real
parameters. See, for instance, the book  \cite[pp. 121--131]{Hale}
for further information about the Hill's equation.
\par We consider the following differential system in $\mathbb{R}^3$:
\begin{equation}
\begin{array}{lll}
\displaystyle \dot{x} & = & - y + z \, x/2,  \vspace{0.2cm} \\
\displaystyle \dot{y} & = & x + z \, y/2,   \vspace{0.2cm} \\
\displaystyle \dot{z} & = & \left(- 2 q (x^2-y^2) - a \right) \, ( x^2+ y^2
-1)\,+z^2,
\end{array} \label{eq6}
\end{equation}
where $a,q \in \mathbb{R}$. This system has the $2 \pi$--periodic
orbit $\Gamma:= \left\{ \gamma(t) \, | \, 0\leq t < 2 \pi
\right\}$ with $\gamma(t) = (\cos t, \sin t, 0)$. \par We consider
the surfaces given by $f_i(x,y,z) =0$, $i=1,2$, with
$\mathbf{f}(x,y,z) = \left( x^2+y^2-1, z \right)^{\rm T}.$ Their
intersection gives the periodic orbit $\Gamma$ and the crossings
over it are transversal, as the following computation shows:
\[ \left| \begin{array}{c} \nabla f_1 (\gamma(t)) \vspace{0.2cm}
\\ \nabla f_2 (\gamma(t)) \vspace{0.2cm} \\ \mathbf{X} (\gamma(t))
\end{array} \right| \, = \, \left| \begin{array}{ccc} 2 \cos(t) &
2 \sin(t) & 0 \vspace{0.2cm} \\ 0 & 0 & 1 \vspace{0.2cm} \\ -
\sin(t) & \cos(t) & 0 \end{array} \right| \, = \, -2 \, \neq \, 0.
\] \par
We have that $ D\mathbf{f}(\mathbf{x}) \, \mathbf{X}(\mathbf{x}) =
\mathbf{k}(\mathbf{x}) \, \mathbf{f}(\mathbf{x}),$ with the following matrix of
cofactors:
\[ \mathbf{k}(x,y,z):= \left( \begin{array}{cc} 0 &  x^2 + y^2 \\ - 2 q (x^2-y^2) - a & z \end{array} \right). \]
Therefore, the cofactor matrix over the periodic orbit reads for:
\[ \mathbf{k}(\gamma(t)) = \left( \begin{array}{cc} 0 & 1 \\
M(t) & 0 \end{array} \right), \] where we denote $M(t):=- 2 q \cos
(2 t) - a$. \par By Theorem \ref{th1}, the stability of $\Gamma$
is given by the eigenvalues of $\mathbf{v}(2 \pi)$, where
$\mathbf{v}(t)$ is the fundamental matrix solution of \[ \left(
\begin{array}{cc} v_{11}'(t) & v_{12}'(t) \\ v_{21}'(t) &
v_{22}'(t) \end{array} \right) \, = \, \left( \begin{array}{cc} 0
& 1 \\ M(t) & 0
\end{array} \right) \, \left( \begin{array}{cc} v_{11}(t) &
v_{12}(t)
\\ v_{21}(t) & v_{22}(t) \end{array} \right). \] Hence,
\begin{itemize}
\item $v_{11}''(t)=M(t) \, v_{11}(t)$ which is the Mathieu's
equation with initial conditions $v_{11}(0)=1$ and $v_{11}'(0)=0$,
and

\item $v_{12}''(t) = M(t) \, v_{12}(t)$ which is the Mathieu's
equation with initial conditions $v_{12}(0)=0$ and $v_{12}'(0)=1$.
\end{itemize}
We have that the system $\mathbf{v}'(t)=\mathbf{k}(\gamma(t)) \,
\mathbf{v}(t)$ coincides with the characteristic system associated
to the Mathieu's equation. Then the stability of the $2
\pi$--periodic orbit $\Gamma:= \{ \gamma(t) \, | \, 0\leq t < 2
\pi \}$ with $\gamma(t) = (\cos t, \sin t, 0)$ of system
(\ref{eq6}) coincides with the stability of the Mathieu's
equation, which is studied in \cite{Simo} and \cite[pp.
128--130]{Hale}.
\newline

{\bf Example 3.} Our third example consists on the study of a
periodic solution related to rigid body dynamics encountered by
Steklov \cite{Steklov}. We consider the motion of a rigid body
around a fixed point in a uniform gravity field. We denote the
weight of the body by $W$ and the distance between the center of
gravity and the fixed point by $\ell$. As described in the work
\cite{Markeev}, we can consider two frames of reference both with
the origin at the moving body. The first frame $OXYZ$ is fixed and
has axis $OZ$ vertical and upward directed. The second frame of
reference is moving solidary with the body and its axes $Ox$, $Oy$
and $Oz$ are directed along the major axes of inertia for the
point $O$, with corresponding moments of inertia denoted by $a$,
$b$ and $c$. We denote by $p$, $q$ and $r$ the components of the
angular velocity vector of the movement and by $\gamma_1$,
$\gamma_2$ and $\gamma_3$, the components of the unit vector in
the direction $OZ$ written in coordinates $Oxyz$. Steklov
considered the case in which the center of mass is located on the
major axis of inertia, which we assume to be the $Ox$--axis. The
Euler--Poisson equations describe the motion of the rigid body:
\begin{equation}
\begin{array}{c} \displaystyle \dot{p} \, = \, \frac{(b-c)}{a} \,
q \, r, \quad \dot{q} \, = \, \frac{(c-a)}{b} \, p \, r \, + \,
\frac{W \ell}{b} \, \gamma_3, \quad \dot{r} \, = \,
\frac{(a-b)}{c} \, p \, q \, - \, \frac{W \ell}{c} \, \gamma_2,
\vspace{0.2cm} \\
\dot{\gamma_1} \, = \, r \, \gamma_2 \, - \, q \, \gamma_3, \qquad
\dot{\gamma_2} \, = \, p \, \gamma_3 \, - \, r \, \gamma_1, \qquad
\dot{\gamma_3} \, = \, q \, \gamma_1 \, - \, p \, \gamma_2.
\end{array}
\label{eq7}
\end{equation}
We are going to assume, without loss of generality and following
\cite{Markeev}, that:
\begin{equation}
b>c, \qquad a+b>c, \qquad b+c>a, \qquad c+a>b, \qquad b>a>2c.
\label{eq8}
\end{equation}
The general Euler--Poisson equations exhibit three first
integrals, whose expressions are described in \cite{bormam}. We
are going to rewrite them for the particular case (\ref{eq7}). The
first one corresponds to the projection of the angular momentum
onto the vertical:
\[ H_1 \, := \, a \, p \, \gamma_1 \, + \, b \, q \, \gamma_2 \, +
\, c \, r \, \gamma_3. \] The second first integral that we encounter is the
geometric property of the vector $(\gamma_1, \gamma_2, \gamma_3)$ to be of
constant modulus:
\[ H_2 \, := \, \gamma_1^2 \, + \, \gamma_2^2 \, + \, \gamma_3^2 \,
- \, 1. \] We also find the full energy (sum of kinetic and
potential energies) of the body as first integral:
\[ H_3 \, := \, \frac{1}{2} \left( a \, p^2 \, + \, b \, q^2 \, +
\, c \, r^2 \right) \, + \, W \ell \, \gamma_1 + \, \frac{(a^2-2ab-2ac+2bc)\,
W\ell}{2(b-a)(a-c)} . \] Straightforward computations show that, if we denote
by $\mathbf{X}$ the vector defined by system (\ref{eq7}), we have $\nabla H_i
\cdot \mathbf{X}^{\rm T} \, \equiv 0$, for $i=1,2,3$. \par As described in
\cite{Markeev}, Steklov in \cite{Steklov} looked for real particular solutions
of system (\ref{eq7}) satisfying the two relations $\gamma_2 \, = \, \beta_2 \,
p \, q$ and $\gamma_3 \, = \, \beta_3 \, r  p$, with $\beta_2$ and $\beta_3$
suitable constants, and he found the following particular periodic solution for
system (\ref{eq7}). Let us define the constant $\mu \, := \, \sqrt{W \ell /
a}$, the dimensionless ``time'' $\nu \, := \, \mu \left( t + t_0 \right)$, with
$t_0$ an arbitrary constant, and the following constants and variable:
\[ \beta_0 \, := \, \frac{a(a-2c)}{(b-a)(a-c)}, \quad \beta_1 \, := \, \frac{a(2b-a)}{(b-a)(a-c)},
\quad  k^2 \, := \, \frac{ b-a }{b-c}, \quad z \, := \,
\frac{1}{k} \, \sqrt{ \frac{a}{a-c}} \, \nu. \] Then, the periodic
orbit found by Steklov can be written as:
\begin{equation}
\begin{array}{ll}
\displaystyle p(t) \,  := \,  \displaystyle - \mu
\sqrt{\frac{\beta_0 (2b-a)}{(a-c)}} \, \, \cn(z;k), &
\displaystyle \gamma_1(t)  :=  \displaystyle 1 \, - \,
\frac{a}{a-c} \, \, \cn^2(z;k), \vspace{0.2cm} \\ \displaystyle
q(t)  := \displaystyle \mu \sqrt{\frac{\beta_0 \, a}{(b-c)}} \, \,
\sn(z;k), & \displaystyle \gamma_2(t)  :=  \displaystyle
\sqrt{\beta_1 k} \, \, \sn(z;k) \, \cn(z;k), \vspace{0.2cm} \\
\displaystyle r(t) := \displaystyle \mu \sqrt{\frac{\beta_1
a}{(a-c)}} \, \, \dn(z;k), & \gamma_3(t) := \displaystyle -
\sqrt{\frac{\beta_0 (b-a)}{(a-c)}} \, \, \cn(z;k)\, \dn(z;k),
\end{array} \label{eq9}
\end{equation}
where $\cn(z;k)$, $\sn(z;k)$ and $\dn(z;k)$ are elliptic Jacobi
functions of variable $z$ and module $k$. \par We recall that the
module $k$ must always satisfy $0 \leq k \leq 1$ and that these
elliptic Jacobi functions are defined as: $\cn(z;k) \, = \,
\cos(\am(z;k))$, $\sn(z;k) \, = \, \sin(\am(z;k))$ and $\dn(z;k)
\, = \, \sqrt{1-k^2\sin^2(\am(z;k))}$, where $\am(z;k)$ is the
Jacobi amplitude defined as the inverse function of the elliptic
integral of first kind ${\rm F}(w;k)$, that is, $\am(z;k)=w$ if
and only if ${\rm F}(w;k)=z$.
\par As stated in \cite{Markeev},
the minimal positive period $T$ of this Steklov solution with
respect to the time $t$ is: \[ T \, = \, 4 \, k \, {\rm K}(k) \,
\sqrt{ \frac{a-c}{W \ell}}, \] where ${\rm K}(k)$ is a complete
elliptic integral of the first kind. See \cite{Abramowitz,
Weisstein} for further information about Jacobi elliptic functions
and integrals. We only recall that:
\[ {\rm F}(z;k) \, := \, \int_{0}^{z} \frac{d
\theta}{\sqrt{1-k^2 \, \sin^2(\theta)}} , \quad {\rm K}(k) \, :=
\, {\rm F} \left( \frac{\pi}{2} ; k \right), \] and
\begin{equation}
\cn^2(z;k) \, + \, \sn^2(z;k) \, = \, 1, \qquad \dn^2(z;k) \, + \,
k^2\, \sn^2(z;k) \, = \, 1. \label{eq10}
\end{equation}
\par
Our goal is to study the characteristic multipliers associated to
the periodic orbit given by Steklov. Since the system has three
functionally independent first integrals, it is well known that
the $6\times6$ monodromy matrix computed from the variational
equations associated to the periodic orbit needs to have at least
4 eigenvalues equal to 1. Several works target to the aim of
obtaining the remaining two characteristic multipliers following
the classical first variational analysis. In \cite{bormam}, a
numerical study of this particular periodic solution is done and a
visualization of it is provided. In the work \cite{Markeev}, the
problem of orbital stability of the Steklov solution is
numerically examined. In \cite{Kucher}, it is proved that if
$a>1$, then this Steklov periodic orbit is orbitally unstable. In
this paper we will see how our method allows to reduce the
computation of these characteristic multipliers to the study of a
linear second order differential equation. \newline

We consider the following five hypersurfaces which contain the
periodic orbit given by Steklov: \begin{eqnarray*}  f_1 & := &
\displaystyle \frac{1}{2} \left( a \, p^2 \, + \, b \, q^2 \, + \,
c \, r^2 \right) \, + \, W \ell \, \gamma_1 \, + \, \frac{(a^2-2ab-2ac+2bc)\, W\ell}{2(b-a)(a-c)} \, = \, 0, \vspace{0.2cm} \\
f_2 & := & \displaystyle \gamma_2 \, - \,
\frac{(a-b)(a-c)}{W\ell(a-2c)} \, p \, q \, = \, 0, \vspace{0.2cm}
\\ f_3 & := &  \displaystyle \gamma_3 \, - \, \frac{(a-b)(a-c)}{W\ell(a-2b)} \,
p \, r \, = \, 0, \vspace{0.2cm} \\ f_4 & := & \displaystyle
\frac{(a-b)}{(a-2c)} \, p^2 \, + \, \frac{(b-c)}{a} \, r^2 \, - \,
\frac{W \ell (a-2b)}{(a-b)(a-c)} \, = \, 0, \vspace{0.2cm} \\ f_5
& := & \displaystyle -\frac{(a-c)}{(a-2b)} \, p^2 \, + \,
\frac{(b-c)}{a} \, q^2 \, + \, \frac{W \ell (a-2c)}{(a-b)(a-c)} \,
= \, 0. \end{eqnarray*} The first hypersurface corresponds to one
of the described first integrals. The second and third
hypersurfaces correspond to the ones looked for by Steklov and the
last two hypersurfaces also contain the periodic orbit and are
independent from the previous ones. Note that only the first one
of the above hypersurfaces is invariant by the flow of
(\ref{eq7}). When we substitute the expressions of the
parameterization (\ref{eq9}), and using (\ref{eq10}), we get that
each $f_i$ values identically zero for any $t$.
\par The expressions of the other two first integrals in relation with the five
polynomials $f_i$ are:
\begin{equation} \begin{array}{l} \displaystyle H_1 \, = \, \displaystyle \frac{a \, p}{ W \ell} \, f_1
\, + \, b \, q \, f_2 \, + \, c \, r \, f_3 \, + \, \frac{ac
\left(a(a-2c) + 2bc\right)p}{2(a-2b)(b-c) W \ell} \, f_4 \, +
\vspace{0.2cm} \\ \displaystyle \quad \displaystyle + \,
\frac{ab \left(a(a-2b) + 2bc\right)p}{2(a-2c)(b-c) W \ell} \, f_5, \vspace{0.3cm} \\
\displaystyle H_2  =  \displaystyle \left( \frac{2 \gamma_1}{W
\ell}  -  \frac{f_1}{W^2 \ell^2} \right)  f_1 \, + \,
\left(\frac{2\delta p \, q}{(a-2c)}  +  f_2 \right) f_2 \, + \,
\left(\frac{2\delta p \, r}{(a-2b)}  +  f_3 \right) f_3 \, +
\vspace{0.2cm}
\\ \quad \displaystyle + \, \left( \frac{-\, a \, c}{(b-c) W\ell}
\, + \, \frac{a\delta \left(a^3-a^2b-2a^2c + abc + 2ac^2 - 2bc^2
\right)
p^2}{(a-2b)^2(a-2c)(b-c)W \ell} \, + \right. \vspace{0.2cm} \\
\qquad \displaystyle \left. \, + \,
\frac{a^2c^2}{4(b-c)^2W^2\ell^2} \, f_4 \, + \,
\frac{ac(4ab-3ac-2bc+2c^2)}{4(b-c)^2W^2\ell^2} \, f_5 \right) \,
f_4 \, + \vspace{0.2cm} \\ \quad \displaystyle + \, \left( \frac{-
\, a \, b}{(b-c) W\ell} \, + \, \frac{a\delta
\left(a^3-2a^2b+2ab^2-a^2c+abc-2b^2
c\right)p^2}{(a-2b)(a-2c)^2(b-c)W \ell} \, + \right. \vspace{0.2cm} \\
\qquad \displaystyle \left. \, + \,
\frac{a^2b^2}{4(b-c)^2W^2\ell^2} \, f_5 \, - \,
\frac{ac(2ab-3ac-2bc+2c^2)}{4(b-c)^2W^2\ell^2} \, f_4 \right) \,
f_5, \end{array} \label{eq11} \end{equation} where $\delta =
(a-b)(a-c)/ (W \ell)$. We note that over the periodic orbit, each
one of these first integrals does not conform a unique leaf since
its coefficients over the hypersurfaces $f_i=0$ are not constants.
This assertion means that the intersection of any combination of
four of the five hypersurfaces $f_i=0$, $i=1,2,3,4,5$, and the set
$H_1=0$ (equivalently $H_2=0$) contains more points than $\Gamma$.
This is the reason why we do not directly use these first
integrals in the computations. We have computed the previous
expressions because they will lead us to show the existence of two
characteristic multipliers equal to $1$.
\newline

The crossings of the five hypersurfaces $f_i=0$, $i=1,2,3,4,5$ are
normal over the periodic orbit because the value of the
determinant of the matrix formed by the gradients of each of the
$f_i$ in the first five rows and the vector field in the last row,
all of them evaluated over the periodic orbit (\ref{eq9}) is: \[
\begin{array}{l} \displaystyle \det \, = \,
\frac{2(b-c)W^2\ell^2}{a(a-c)^2} \left(
\frac{2a(a-2c)(a-c)}{(b-a)^2} \, + \,
\frac{2a(a-c)^2}{(2b-a)(b-a)} \, p(t)^2 \, + \right. \vspace{0.2cm} \\
\displaystyle \left. \quad + \, \frac{2(2b-a)(a-c)(b-c)}{a(b-a)}
\, q(t)^2 \,  + \,
\frac{2(2b+2c-3a)(a-c)^2(b-c)}{(2b-a)(b-a)W\ell} \, r(t)^2 \, +
\right. \vspace{0.2cm} \\ \displaystyle \left. \quad
 + \, \frac{8(a-c)(b-c)W\ell}{a(2b-a)} \, \gamma_2(t)^2 \right),
\end{array}
\] which is positive for any $t$. This assertion is true because the conditions (\ref{eq8}) imply that
all the coefficients are positive except the term $(2b+2c-3a)$,
which can be positive, negative or zero. If it is positive or
zero, we already have that $\det >0$. If it is negative, we are
going to show that: \[ \frac{2a(a-2c)(a-c)}{(b-a)^2} \, - \,
\frac{2(3a-2b-2c)(a-c)^2(b-c)}{(2b-a)(b-a)W\ell} \, r(t)^2 \, > \,
0 \] which ensures that $\det
>0$ for any value of $t$. We note that $-1 \leq \dn(z;k) \leq 1$, so we consider
any value of $t$ for which $\dn(z;k)^2$ is equal to $1$ and, using
(\ref{eq9}) and some computations, we can bound the previous
expression by:
\[ \begin{array}{l} \displaystyle \frac{2a(a-2c)(a-c)}{(b-a)^2} \, - \,
\frac{2(3a-2b-2c)(a-c)^2(b-c)}{(2b-a)(b-a)W\ell} \, r(t)^2 \, \geq
\, \vspace{0.2cm} \\ \displaystyle \quad \, \geq \,
\frac{2a(a-2c)(a-c)}{(b-a)^2} \, - \,
\frac{2(3a-2b-2c)(a-c)^2(b-c)}{(2b-a)(b-a)W\ell} \, \frac{W \ell
\, (2b-a)a}{(b-a)(a-c)^2} \, = \, \vspace{0.2cm} \\
\displaystyle \quad \quad \, = \, \frac{2a\, (2b-a)}{(b-a)} \, >
\, 0.
\end{array} \]
\par
The matrix of cofactors $\mathbf{k} = \left(k_{ij}\right)$
associated to the five previous hypersurfaces, that is,
\[ D\mathbf{f}\,\mathbf{X}\,=\,\left( \begin{array}{l} \nabla f_1 \\ \nabla f_2 \\ \nabla f_3
\\ \nabla f_4 \\ \nabla f_5 \end{array} \right) \, \mathbf{X}
\, = \, \mathbf{k} \, \left( \begin{array}{c} f_1 \\
f_2 \\ f_3 \\ f_4 \\ f_5 \end{array} \right), \] reads for:
\[ \begin{array}{l}
\displaystyle k_{1j} \, := \, 0 \ \quad \mbox{for} \quad \
i=1,2,\ldots ,5, \vspace{0.3cm} \\
\displaystyle k_{21} \, := \, \frac{-\, r}{W\ell}, \quad k_{22} \,
:= \, 0, \quad k_{23} \, := \, -\frac{(a^2-2ab-ac+3bc)\, p}{b(a-2c)}, \vspace{0.2cm} \\
\displaystyle  k_{24} \, := \, \frac{a c \, r}{2(b-c)W\ell}, \ \,
k_{25} \, := \, - \frac{(a^2b - 2ab^2 - 2a^2c + 2abc + 2b^2c +
2ac^2 - 2bc^2)\, r}{2(a-2c)(b-c) W \ell}, \vspace{0.3cm} \\
\displaystyle k_{31} \, := \, \frac{q}{W \ell}, \quad k_{32} \, :=
\, \frac{(a^2-ab-2ac+3bc) \, p}{(a-2b)c}, \quad k_{33} \, := \, 0, \vspace{0.2cm} \\
\displaystyle  k_{34} \, := \,- \frac{(2a^2b - 2ab^2 - a^2c - 2abc
+ 2b^2c + 2ac^2-2bc^2)\, q}{2(a-2b)(b-c) W \ell}, \ \, k_{35}
\, := \, \frac{-ab \, q}{2(b-c)W\ell}, \vspace{0.3cm} \\
\displaystyle k_{41} \, := \, k_{43} \, := \,  k_{44} \, := \,
k_{45} \, := \, 0, \quad k_{42} \, := \, \frac{-2(b-c)W\ell \,
r}{ac}, \vspace{0.2cm} \\
\displaystyle k_{51} \, := \, k_{52} \, := \,  k_{54} \, := \,
k_{55} \, := \, 0, \quad k_{53} \, := \, \frac{2(b-c)W\ell \,
q}{ab}.
\end{array} \] The fundamental matrix of solutions,
evaluated in $T$, of the linear differential system of equations
\begin{equation} \mathbf{v}'(t) \, = \, \mathbf{k}(\gamma(t))
\, \mathbf{v}(t), \label{eq12}
\end{equation} with, \[ k \left( \gamma(t) \right) \, =  \, \left(
\begin{array}{ccccc} 0 & 0 & 0 & 0 & 0  \\ k_{21}(t) & 0 & k_{23}(t) & k_{24}(t) & k_{25}(t)
\\ k_{31}(t) & k_{32}(t) & 0 & k_{34}(t) & k_{35}(t) \\ 0  & k_{42}(t) & 0 & 0 & 0
\\ 0 & 0 & k_{53}(t) & 0 & 0 \end{array} \right) \]
the matrix of cofactors evaluated over the periodic orbit
(\ref{eq9}), gives the characteristic multipliers associated to
the periodic orbit encountered by Steklov. We note that the first
equation of this system is $v_1'(t)=0$ which corresponds to the
fact that $f_1$ is a first integral of the system and gives us
that $v_1(t)$ needs to be a constant. Analogously, the other two
first integrals give us relations among the values of the
solutions of system (\ref{eq12}). In particular, from
(\ref{eq11}), $H_1$ gives us that the function:
\[ \begin{array}{l} \displaystyle h_1(t) \, := \, \displaystyle \frac{a \, p(t)}{ W \ell} \,
v_1(t) \, + \, b \, q(t) \, v_2(t) \, + \, c \, r(t) \, v_3(t) \,
+  \vspace{0.2cm} \\ \displaystyle \quad \displaystyle + \,
\frac{ac \left(a(a-2c) + 2bc\right)p(t)}{2(a-2b)(b-c) W \ell} \,
v_4(t) \, + \, \frac{ab \left(a(a-2b) +
2bc\right)p(t)}{2(a-2c)(b-c) W \ell} \, v_5(t), \end{array} \]
satisfies $h_1'(t) \equiv 0$ when considered over the solutions of
system (\ref{eq12}). And the first integral $H_2$ gives the
function:
\[ \begin{array}{l} h_2(t) \, := \,  \displaystyle \left( \frac{2 \gamma_1(t)}{W
\ell} \right) \, v_1(t) \, + \, \left(\frac{2\delta\,  p(t) \,
q(t)}{(a-2c)} \right) \, v_2(t) \, + \, \left(\frac{2\delta\, p(t)
\, r(t)}{(a-2b)} \right) \, v_3(t) \, +  \vspace{0.2cm} \\
\ \displaystyle +  \, \left( \frac{-\, a c}{(b-c) W\ell} \, + \,
\frac{a\delta \left(a^3-a^2b-2a^2c + abc + 2ac^2 - 2bc^2 \right)
p(t)^2}{(a-2b)^2(a-2c)(b-c)W\ell} \right) v_4(t) \, +
\vspace{0.2cm} \\ \ \displaystyle + \, \left( \frac{- \, a
b}{(b-c) W\ell} \, + \, \frac{a\delta
\left(a^3-2a^2b+2ab^2-a^2c+abc-2b^2
c\right)p(t)^2}{(a-2b)(a-2c)^2(b-c)W\ell} \right) \, v_5(t),
\end{array} \] where $\delta = (a-b)(a-c)/(W \ell)$, which is also constant over the solutions of
(\ref{eq12}). \par We obtain these constant functions over the
solutions of (\ref{eq12}) because each function $H_i(\gamma(t))$,
$i=1,2$, is $0$ for any $t$. We can then derive $H_i(\gamma(t))$
with respect to $t$ and we deduce the functions $h_i(t)$ taking
into account that $v_i(t)$ is related with $\nabla
f_i(\gamma(t))$, $i=1,2,\ldots, 5$.
\newline

We first note that the following constants conform a solution of
system (\ref{eq12}):
\[ \begin{array}{c} \displaystyle v_1(t) \, = \, 1, \quad v_2(t) \, = \, 0,
\quad v_3(t) \, = \, 0, \vspace{0.2cm} \\ \displaystyle v_4(t) \,
= \, \frac{2(a-2b)}{a^2-2ab-2ac+2bc}, \quad v_5(t) \, = \,
\frac{-2(a-2c)}{a^2-2ab-2ac+2bc} . \end{array} \] The functions
$h_1(t)$ and $h_2(t)$ take the constant values $0$ and $
4\delta/(a^2-2ab-2ac+2bc)$, with $\delta = (a-b)(a-c)/(W \ell)$,
over this solution, respectively. Since this solution takes a
nonzero value on $v_1(t)$, which needs to be constant, we may get
the first column of the fundamental matrix solution of
(\ref{eq12}) from the last four columns and this particular
constant solution. Hence, we only need to take into account the
last four columns of the fundamental matrix solution of
(\ref{eq12}). \newline

Let us now consider the fundamental matrix $\mathbf{v}(t)$
solution of (\ref{eq12}) whose initial condition is the identity
matrix. We define $v_{ij}(t)$ as the function corresponding to row
$i$ column $j$ of this fundamental matrix of solutions. Using that
the three functions $v_1(t)$, $h_1(t)$ and $h_2(t)$ are constants,
we can obtain certain relations among the rows of this matrix at
the value $T$. \par We define the constants: \[ \rho_1 \, := \,
\frac{2 (2b-a)(b-c) W \ell}{\sqrt{a} (b-a)(a-2c)^{3/2}}, \quad
\rho_2 \, := \, \frac{2 \sqrt{a} (b-c) W \ell}{(b-a) \sqrt{a-2c}
(a^2-2ab+2bc)}, \] and we get that $\mathbf{v}(T)$ equals
\[ \left( \begin{array}{ccccc} \displaystyle 1 & 0 & 0 & 0
& 0 \vspace{0.2cm} \\ \displaystyle  v_{21}(T) & v_{22}(T) &
v_{23}(T) & v_{24}(T) & v_{25}(T) \vspace{0.2cm} \\ \displaystyle
v_{31}(T) & v_{32}(T) & v_{33}(T) & \displaystyle
\frac{v_{44}(T)-1}{\rho_1} &
\displaystyle \frac{v_{55}(T)-1}{\rho_2} \vspace{0.2cm} \\
\displaystyle \rho_1 v_{31}(T) & \rho_1 v_{32}(T) & \rho_1 \left(
v_{33}(T) -1 \right) & v_{44}(T) & \displaystyle
\rho_1 \left( \frac{v_{55}(T)-1}{\rho_2} \right) \vspace{0.2cm} \\
\displaystyle \rho_2 v_{31}(T) & \rho_2 v_{32}(T) & \rho_2 \left(
v_{33}(T) -1 \right) & \displaystyle \rho_2 \left(
\frac{v_{44}(T)-1}{\rho_1} \right) & v_{55}(T)
\end{array} \right). \] The computation of the characteristic
polynomial associated to this matrix gives:
\[ \det \left(v(T) - \lambda \, I\right) \, = \, -(\lambda-1)^3 \left( \lambda^2 \, + \,
B \, \lambda \, + \, C \right), \] where \[ B \, := \, 2 \, - \,
v_{22}(T) \, - \, v_{33}(T) \, - \, v_{44}(T) \, - \, v_{55}(T),
\] and \begin{eqnarray*} C & := & v_{22}(T) \left( v_{33}(T) \, + \, v_{44}(T) \,
+ \, v_{55}(T) \, - \, 2 \right) \, + \\ & & \ - \, v_{32}(T)
\left( v_{23}(T) \, + \, \rho_1 \, v_{24}(T) \, + \, \rho_2 \,
v_{25}(T) \right). \end{eqnarray*} We note that each one of the
first integrals gives an eigenvalue equal to one. Moreover, using
Liouville's formula, we know that the product of all the
eigenvalues of $\mathbf{v}(T)$ is equal to one because the trace
of matrix $\mathbf{k}(\gamma(t))$ is identically zero and
$\mathbf{v}(t)$ is a fundamental solution of (\ref{eq12}). Hence,
we get that $C=1$. Then, we have that the roots of $\lambda^2 \, +
\, B \,  \lambda \, + \, 1 \, = \, 0$ are:
\begin{itemize}
\item if $B^2 \, > \, 4$, then the roots are both real, and one of
modulus greater than $1$ and the other with modulus lower than
$1$: so the Steklov solution is unstable, \item if $B^2 \, \leq \,
4$, then the roots are both of modulus equal to one and the
characteristic multipliers do not decide if the Steklov solution
is unstable or not. \end{itemize}

Let us go back to the differential linear system (\ref{eq12}). We
know that $v_1(t)$ is constant, and we denote its value by $v_{10}
\, := \, v_{1}(t)$. Analogously, we denote by $h_{10} \, :=
h_1(t)$ and $h_{20} \, := h_2(t)$, the values of the other two
constants. We can determine $v_2(t)$ and $v_3(t)$ from the last
two equations of (\ref{eq12}) and also from the two relations
$h_1(t) \, = \, h_{10}$ and $h_2(t) \, = h_{20}$. We equate the
two expressions of $v_2(t)$ and $v_3(t)$ obtained in these two
different ways and we get a linear differential system of
equations for $v_4(t)$ and $v_5(t)$ which reads for:
\[ \begin{array}{l} \displaystyle \frac{\partial \,
\dn(z;k)}{\partial t} \, v_4'(t) \, - \,
\frac{W\ell(2b-a)\dn(z;k)}{a(a-c)(b-c)} \left[ \left( c -
\frac{ab-2ac+2bc}{2b-a} \, \cn(z;k)^2 \right) v_4(t) \, + \right.
\vspace{0.2cm} \\ \displaystyle \left. \quad + \, \left( b -
\frac{(2b-a)(ab-ac-bc)}{(a-2c)c}\, \cn(z;k)^2 \right) v_5(t) \, +
\right. \vspace{0.2cm} \\ \displaystyle \left. \quad + \,
\frac{2(b-c) \sqrt{W \ell (b-a)(a-2c)}}{ac \sqrt{2b-a}} \,
\cn(z;k) \, h_{10} \, + \, \frac{2(b-c)}{c} \, v_{10} \,
\cn(z;k)^2 \, + \right. \vspace{0.2cm} \\ \displaystyle \left.
\quad + \, \frac{b-c}{a} \left(W \ell \, h_{20} \, - \, 2 v_{10}
\right) \right] \, = \, 0, \vspace{0.3cm} \\ \displaystyle
\frac{\partial \, \sn(z;k)}{\partial t} \, v_5'(t) \, + \,
\frac{W\ell(a-2c)\sn(z;k)}{a(b-a)(b-c)} \left[ \left( c -
\frac{r_0(ab-ac+bc)}{(2b-a)} \cn(z;k)^2 \right) v_4(t) \, +
\right. \vspace{0.2cm} \\ \displaystyle \left. \quad + \, \left( b
+ \frac{(b-a)(2ab-ac-2bc)}{(a-2c)(a-c)} \, \cn(z;k)^2 \right)
v_5(t) \, + \right. \vspace{0.2cm} \\ \displaystyle \left. \quad +
\, \frac{-2(b-c)\sqrt{W \ell (b-a)(2b-a)}}{ab\sqrt{a-2c}} \,
\cn(z;k) \, h_{10} \, - \, \frac{2(b-a)(b-c)}{b(a-c)} \,
\cn(z;k)^2 \, v_{10} \, + \right. \vspace{0.2cm} \\ \displaystyle
\left. \quad + \, \frac{b-c}{a} \left(W \ell \, h_{20} \, - \, 2
v_{10} \right) \right] \, = \, 0, \end{array} \] where $r_0 \, =
\, (b-a)(a-2c)/(b(a-c))$. Equating $v_5(t)$ from the first
equation and substituting its value in the second one, we get a
second order linear differential equation for $v_4(t)$ whose
fundamental set of solutions would let us to the computation of
the characteristic multipliers associated to the Steklov periodic
orbit. We define
\[ \omega(t) \, := \, b -
\frac{(2b-a)\delta_1}{(a-2c)c}\, \cn(z;k)^2,\quad \delta_1 \, :=
\, ab-ac-bc \] and this second order linear differential equation
for $v_4(t)$ is:
\begin{equation}
v_4''(t) \, - \left( \frac{
\omega'(t)}{\omega(t)} \, + \, \frac{2}{\dn(z;k)} \,
\frac{\partial \dn(z;k)}{\partial t} \right) v_4'(t) \, + \,
A_0(t) \, v_4(t) \, + \, A_{nh}(t) \, = \, 0, \label{eq13}
\end{equation}
where
\[ \begin{array}{l} \displaystyle \frac{A_0(t)}{W\ell} \, := \,
\frac{2(b-c)(a^2-4b(a-c))}{a(a-2c)(a-c)}\,
\frac{\dn(z;k)^2}{\omega(t)} \, + \, \frac{4}{a} \,+\,
\frac{a^2-2(\delta_1+2ac)}{bc(a-c)}\, \cn(z;k)^2 , \vspace{0.3cm}
\\ \displaystyle A_{nh}(t) \, := \, \frac{-2 (2b-a)^2}{ a\delta(a-2c)c^2\, \omega(t)} \, \Big( \alpha_3(t) \, h_{10} \, + \,
\alpha_4(t) \, v_{10} \, +
 \, \alpha_5(t)  \left(2
v_{10} - W \ell \, h_{20} \right) \Big),
\end{array}
\] with $\delta\, :=\, (a-b)(a-c)/(W\ell)$ and
\begin{eqnarray*}
\alpha_3(t) & := & \sqrt{\frac{W \ell(b-a)}{(2b-a)(a-2c)}} \left(
\frac{(a-2c)(3ac(b+c-a)+5bc(b-c))}{(2b-a)} \, + \right. \vspace{0.2cm} \\
& & \left. + 2ab(2c-a) \, + \, \frac{\delta_1
\left((b-a)^2+\delta_1\right)}{b} \, \cn(z;k)^2\right) \cn(z;k),
\vspace{0.3cm} \\ \alpha_4(t) & := & -a(a-2c) \left(
\frac{bc}{2b-a} \, + \, (b-2c) \, \cn(z;k)^2 \, + \, \frac{(b-a)\,
\delta_1}{b(a-2c)} \, \cn(z;k)^4 \right),
 \vspace{0.3cm} \\ \alpha_5(t) & := & \frac{a b c}{2} +
 \frac{c^2(a^2-a b-b^2-ac+bc)}{2b-a} \, - \, \frac{(a-2b+2c)
 \, \delta_1}{2} \, \cn(z;k)^2 .
\end{eqnarray*} \newline

In short, we have reduced the problem of computing the
characteristic multipliers for the periodic orbit (\ref{eq9}) to
the study of the second--order linear differential equation
(\ref{eq13}). The classical approach to this problem is to compute
the characteristic multipliers via the first variational
equations, thus involving a linear differential system of order
$6$. Our method starts with   a linear differential system of
order $5,$ see (\ref{eq12}). In this particular problem, we know
three first integrals which let us reduce the order by $3$,
getting the second--order linear differential equation
(\ref{eq13}). We note that the use of the first integrals is not
trivial, since we need to write them in terms of the considered
hypersurfaces, see (\ref{eq11}), and then relate them to the
variables of the linear differential system (\ref{eq12}).


\begin{thebibliography}{99}

\bibitem{Abramowitz} {\sc M. Abramowitz and I.A. Stegun},{\it \ Handbook
of mathematical functions with formulas, graphs, and mathematical
tables.} Reprint of the 1972 edition. Dover Publications, Inc.,
New York, 1992.

\bibitem{bormam} {\sc A. V. Borisov and I. S. Mamaev}, {\it \ Euler-Poisson equations and integrable cases},
Regul. Chaotic Dyn. {\bf 6} (2001), 253--276.

\bibitem{Simo} {\sc H. Broer and C. Sim\'o}, {\it \ Resonance tongues in Hill's equations: a
geometric approach},  J. Differential Equations {\bf 166} (2000),
290--327.

\bibitem{dekleine} {\sc H. A. De Kleine}, {\it \ A note on the asymptotic stability of
periodic solutions of autonomous differential equations.} SIAM
Rev. {\bf 26} (1984), 417--421.

\bibitem{hyper} {\sc H. Giacomini and M. Grau}, {\it \ On the stability of
limit cycles for planar differential systems}, J. Differential
Equations {\bf 213} (2005), 368--388.

\bibitem{Hale} {\sc J. K. Hale}, {\it \ Ordinary differential equations.} Pure and Applied Mathematics,
Vol. {\bf XXI}. Wiley--Interscience (John Wiley \& Sons), New
York--London--Sydney, 1969.

\bibitem{Kucher} {\sc E. Yu. Kucher}, {\it \ Characteristic exponents of periodic Steklov and Chaplygin
solutions}. (Russian) Mekh. Tverd. Tela  {\bf 33} (2003), 33--39.

\bibitem{Liapunov} {\sc A. M. Liapunov}, {\it \ Stability of motion.} Translated from the Russian
by F. Abramovici and M. Shimshoni. Mathematics in Science and
Engineering, Vol. {\bf 30}. Academic Press, New York--London,
1966.

\bibitem{Markeev} {\sc A. P. Markeev}, {\it \ On the Steklov case in rigid body
dynamics}, Regul. Chaotic Dyn. {\bf 10} (2005), 81--93.

\bibitem{Perko} {\sc L. Perko},{\it \ Differential equations
and dynamical systems.} Third edition. Texts in Applied
Mathematics, {\bf 7}. Springer--Verlag, New York, 2001.

\bibitem{Steklov} {\sc V. A. Steklov}, {\it \ New particular
solution of differential equations of motion of a heavy rigid body
about a fixed point}, Trudy Ob-va estest. (1899), {\bf 1}, 1--3.

\bibitem{Weisstein} {\sc E. W. Weisstein}, {\it \ Jacobi Elliptic Functions.} From
MathWorld--A Wolfram Web Resource. \\ \texttt{
http://mathworld.wolfram.com/JacobiEllipticFunctions.html}

\end{thebibliography}
\end{document}